\documentclass[10pt,letterpaper,conference]{ieeeconf}

\usepackage{latexsym, amssymb, amsmath}
\usepackage[dvips]{graphicx} % added for inserting the pictures

\usepackage{float}      %these two packages are for fixing figures in place
\restylefloat{figure}
\usepackage{graphicx,dsfont}
\usepackage{mathrsfs}  % package for math fonts: mathscripts \mathscr
\usepackage{color} % standard color package
\usepackage[latin1]{inputenc}
\usepackage{enumerate}
\usepackage{flushend} % align at the last lines of every page
\usepackage{stfloats}
\usepackage{cite} % package to reduce the length of long citations by using a hyphen

\usepackage{balance}
\allowdisplaybreaks[4] % added by Hong for allowing page breaks among a large formula

\newtheorem{prop}{Proposition} % [section]
\newtheorem{proposition}[prop]{Proposition} % [section]
\newtheorem{lem}{Lemma}%[section]
\newtheorem{lemma}[lem]{Lemma} % [section]
\newtheorem{thm}{Theorem} % [section]
\newtheorem{theorem}[thm]{Theorem} % [section]
\newtheorem{cor}{Corollary} % [section]
\newtheorem{corollary}[cor]{Corollary} % [section]
\newtheorem{defn}{Definition} % [section]
\newtheorem{definition}[defn]{Definition} % [section]
\newtheorem{exmp}{Example} % [section]
\newtheorem{example}[exmp]{Example} % [section]
\newtheorem{exam*}{Example}

\def\custombibliography#1{
 \normalsize
\section*{\centering References}
 \list
 {[\arabic{enumi}]}{\settowidth\labelwidth{[#1]}\leftmargin\labelwidth
 \setlength{\itemsep}{.1em}
 \advance\leftmargin\labelsep
 \usecounter{enumi}}
 \def\newblock{\hskip .11em plus .33em minus -.07em}
 \sloppy
 \sfcode`\.=1000\relax}

\def\L2{{\cal L}_2}
\def\bull{\rule{0.08in}{0.08in}} % square filled bullet
\def\openbull{\framebox[0.08in][c]{$\;$}} % square unfilled bullet
\def\re{{\mathbb R}} % real numbers (AMS symbol)

\def\C{{\mathbb C}} % complex numbers (AMS symbol)
\def\ss#1{{\scriptstyle #1}} % fonts in subsubscripts
\def\eqref#1{(\ref{#1})} % parentheses around referenced equation numbers
\def\Tr{{\rm Tr}}
\def\vec{{\rm vec}}

\newcommand{\comment}[1]{} % Allows one to comment out a block of text

\def\begce{\begin{center}}
\def\endce{\end{center}}
\def\begar{\begin{array}}
\def\endar{\end{array}}
\def\begeq{\begin{equation}}
\def\endeq{\end{equation}}
\def\begdi{\begin{displaymath}}
\def\enddi{\end{displaymath}}
\def\begdis{\begin{eqnarray*}}
\def\enddis{\end{eqnarray*}}
\def\begeqa{\begin{eqnarray}}
\def\endeqa{\end{eqnarray}}
\def\begdes{\begin{description}}
\def\enddes{\end{description}}
\def\begit{\begin{itemize}}
\def\endit{\end{itemize}}
\def\begen{\begin{enumerate}}
\def\enden{\end{enumerate}}
\def\beglar{\left[\begin{array}}
\def\endrar{\end{array}\right]}
\def\begle{\begin{lemma}}
\def\endle{\end{lemma}}
\def\begde{\begin{definition}}
\def\endde{\end{definition}}
\def\begth{\begin{theorem}}
\def\endth{\end{theorem}}
\def\begco{\begin{corollary}}
\def\endco{\end{corollary}}
\def\begprop{\begin{proposition}}
\def\endprop{\end{proposition}}
\def\begex{\begin{example}}
\def\endex{\hfill\openbull \end{example} \vspace*{0.1in}}
\def\begexer{\begin{exercise}}
\def\endexer{\end{exercise}}

\def\begres{\noindent{\bf Remarks}:\begin{enumerate}}
\def\endres{\end{enumerate} \par}
\def\begpr{\noindent{\em Proof:}$\;\;$}
\def\endpr{\hfill\bull \vspace*{0.1in}}

\def\begtab{\begin{tabular}}
\def\endtab{\end{tabular}}
\def\rref#1{(\ref{#1})}
\newcommand\cdcout[1]{} % reduce it to 6 pages
\newcommand{\rv}[1]{\boldsymbol{#1}} % Use italic boldface to indicate the Random Variables
\newcommand{\RomanNumber}[1]{\uppercase\expandafter{\romannumeral #1}}
\newcommand{\romannumber}[1]{\lowercase\expandafter{\romannumeral #1}}

\DeclareMathAlphabet{\mathpzc}{OT1}{pzc}{m}{it}

\def\1{\rv 1} %Indicator
\def\allseriesA{\mbox{$\re\langle\langle A\rangle\rangle$}}
\def\allseriesA'{\mbox{$\re\langle\langle A'\rangle\rangle$}}

\def\L1spaceprodu{{ L}_1(\Omega\times [0,T],{\mathcal P},P\otimes \lambda)}

\def\Hspace0{{\mathcal H}^2_0}

\def\ket#1{\left| #1 \right\rangle}
\def\bra#1{\left\langle #1 \right|}

\usepackage{cite} % generates citation ranges
\usepackage{color} % standard color package
\usepackage{graphicx} % standard graphics package
\usepackage{latexsym,amssymb,amsmath} % extended symbol packages
\usepackage{mathrsfs} % mathscripts \mathscr
\usepackage{multicol}
\usepackage[german,english]{babel}
\usepackage[T1]{fontenc}
\usepackage[latin1]{inputenc}

\allowdisplaybreaks[4]

\IEEEoverridecommandlockouts                              % This command is only
\title{Physical Realizability of Multi-Level Quantum Systems$^\ast$\thanks{$^\ast$This research was supported by the Australian Research Council.}}

\author{Luis~A.~Duffaut~Espinosa$^\dagger$\thanks{$^\dagger$School of Engineering and Information Technology, University of New South Wales at ADFA, Canberra, ACT 2600, Australia. {\tt\small \{l.duffaut,  i.petersen, v.ougrinovski\}@adfa.edu.au.}}, Z.~Miao$^\ddag$\thanks{$^\ddag$Research School of Information Sciences and Engineering, Canberra, ACT 2601, Australia. {\tt\small zibo.miao@anu.edu.au}.}, I.~R.~Petersen$^\dagger$, V.~Ugrinovskii$\,^\dagger$, and M.~R.~James$^\S$\thanks{$^\S$ARC Centre for Quantum Computation and Communication Technology, Research School of Engineering, Australian National University, Canberra, ACT 0200, Australia. {\tt\small matthew.james@anu.edu.au}.}}

\date{\today}

\begin{document}

\maketitle

\begin{abstract}
This paper considers the physical realizability condition for multi-level quantum systems having polynomial Hamiltonian and multiplicative coupling with respect to several interacting boson fields. Specifically, it generalizes a recent result the authors developed for two-level quantum systems. For this purpose, the algebra of $\pmb{SU(n)}$ was incorporated. As a consequence, the obtained condition is given in terms of the structure constants of $\pmb{SU(n)}$.
\end{abstract}

\section{Introduction}  \label{sec:section1}

In an environment where the classical laws of physics apply, standard control techniques such as optimization or a Lyapunov procedures do not worry in general of the nature of the controller they synthesized. In other words, their implementation is always possible since the physics behind them still hold. However, if one desires to implement a controller that obeys the laws imposed by quantum mechanics (e.g., quantum coherent control \cite{Sarovar-Ahn-Jacobs-Milburn_2004,Helon-James_2006,Lloyd_2000}), then such a task is not so easily achieved unless an explicit characterization of those laws is given in terms of the control system vector fields. %Moreover, if a quantum system is controlled by a classical controller, then undesired problems such as acquisition of suitable quantum information, quantum error correction, etc arise \cite{Sarovar-Ahn-Jacobs-Milburn_2004,Helon-James_2006,Lloyd_2000}. 
This is exactly the purpose for introducing the concept of \emph{physical realizability}. 

Conditions for physical realizability were first given specifically for linear systems satisfying the quantum harmonic oscillator canonical commutation relations \cite{James-Nurdin-Petersen_2008,Maalouf-Petersen_2011}. Recently, the formalism was extended for systems describing the dynamics of open two-level quantum systems interacting only with one quantum field in which the algebra of $SU(2)$ played a central role \cite{Duffaut-et-al_2012b}. Compared to a linear quantum system, the systems being analyzed were bilinear, and the the commutation relations were dependent on the system variables. Thus, the main contribution of this paper, given in Section \ref{sec:section4}, is to provide a condition for physical realizability of multi-level quantum systems having polynomial Hamiltonian and multiplicative coupling, and whose system variables obey the commutation relations described by the algebra of $SU(n)$. As expected, the obtained condition is given in terms of the symmetric and antisymmetric structure 
constants of $SU(n)$. Another contribution is that the systems under consideration have been allowed to interact with multiple quantum fields in quadrature form.

The paper is organized as follows. Section \ref{sec:section2} presents the basic preliminaries on open quantum systems. In Section \ref{sec:section3}, the necessary algebraic machinery to study open multi-level quantum systems is given. This is followed by Section \ref{sec:section4}, in which the definition of physical realizability is provided as well as a condition for a bilinear QSDE to be physically realizable. Finally, Section \ref{sec:conclusions} gives the conclusions.

\section{Open Multi-Level Quantum Systems} \label{sec:section2}

\emph{Open quantum systems} are systems governed by the laws of quantum mechanics that interact with an external environment. % (e.g., electromagnetic field) are known as \emph{open quantum systems}. 
%In order to study such systems, one has to give a quantum description of both the system and the interacting environment. 
A quantum mechanical system is described in terms of \emph{observables} and \emph{states}. Observables represent physical quantities that can be measured, as self-adjoint operators on a complex separable Hilbert space $\mathfrak{H}$, while states give the current status of the system, as elements of $\mathfrak{H}$, allowing the computation of expected values of observables. In \cite{Bouten-Handel-James_2007,Parthasarathy_92}, the evolution of open quantum systems is given in terms of quantum stochastic differential equations. For this purpose, observables may be thought as quantum random variables that do not in general commute. A measure of the non commutativity between observables is usually given by the \emph{commutator} between operators. 
\begde The commutator of two scalar operators $x$ and $y$ in ${\mathfrak{H}}$ is the antisymmetric bilinear operation  
\begdi    
[x, y] = xy - yx.
\enddi
Also, if $x$ is an $n_1$-dimensional vector of operators in ${\mathfrak H}$ 
% , the commutator of ${x}$ and a scalar operator $y$ in ${\mathfrak{H}}$ is the $n_1$-dimensional vector of operators %in ${\cal H}$ 
% %\begdi
% $[{x},y] \triangleq {x} y - y {x}$,  
% %\enddi 
% and the commutator of ${x}$ 
and $y$ is an $n_2$-dimensional vector of operators in ${\mathfrak H}$, then % is the $n_1\times n_2$ matrix of operators 
\begdi
[{x},{y}^T] \triangleq {x} {y}^T - ({y} {x}^T)^T,
\enddi %\vspace*{-0.2in}
which is an $n_1\times n_2$ matrix of operators in ${\mathfrak H}$.
\endde

\noindent This commutator satisfies
\begeq \label{eq:commutator_vector_transpose}
[x,y^T]^T = -yx^T + \left(x y^T \right)^T = -[y,x^T]. 
\endeq  
The adjoint of $x$ is denoted by $x^\dagger = (x^\#)^T$ with \begdi {x}^\# \triangleq \left( \begin{array}{c} x_1^\ast \\ x_2^\ast \\ \vdots \\ x_n^\ast \end{array}\right) \enddi and $^\ast$ denotes the operator adjoint. % is the adjoint of the $n$-dimensional complex vector ${x}=(x_1\; x_2 \;\cdots\; x_n )^T$. 
In the case of complex vectors and matrices, $^\ast$ denotes the complex conjugate while $^\dagger$ denotes the conjugate transpose. The non-commutativity of observables is a fundamental difference between quantum systems and classical systems in which the former must satisfy certain commutation relations originating from \emph{Heisenberg uncertainty principle}. 

The environment consists of a collection of oscillator systems each with annihilation field operator $w(t)$ and creation field operator $w^\ast (t)$ used for the annihilation and creation of quanta at point $t$, and commonly known as the \emph{boson quantum field} (with parameter $t$). Here it is assumed that $t$ is a real time parameter. The field operators $w(t)$ and $w^\ast(t)$ satisfy commutation relations as well. That is, \begdi [w(t), w^\ast (t')] = \delta(t - t' ) \enddi for all $t, t'\in \re$, where $\delta(t)$ denotes the Dirac delta. Its mathematical description is given in terms of a Hilbert space called a \emph{Fock space}.~When the boson quantum field is in the vacuum state, i.e., no physical particles are present, it then represents a natural quantum extension of white noise, and may be described using the quantum It\^o calculus \cite{Bouten-Handel-James_2007,Parthasarathy_92}.~This amounts to having three interacting signals (inputs) in the evolution of the system: the annihilation processes $W(t)$, the 
creation process $W^\ast (t)$, and the counting process $\Lambda_w(t)$. 

In simple words, the evolution of an open quantum system is described putting together the evolutions of the system and the environment in an unitary fashion. That is, if $\psi$ is an initial state then 
%\begdi
$\psi(t) =U(t) \psi$, %and %\;\; 
%$x(t)=U^\dagger(t) \psi U(t)$, 
%\enddi
where $U(t)$ is unitary for all $t$, and is the solution of 
\begin{align*}
dU(t)  = &\left(\rule{0in}{0.18in} (S - I)\,d\Lambda_w(t) + L \,d W^\ast(t) - L^\ast S\,d W(t)\right. \\ 
& {}\;\;\; \left. - \frac{1}{2}( L^\ast L + \pmb{i} \mathcal{H})\,dt \right) U(t),
\end{align*}
with initial condition $U(0) = I$, $I$ denoting the identity operator and $\pmb{i}$ being the imaginary unit. Here, $\mathcal{H}$ is a fixed self-adjoint operator representing the \emph{Hamiltonian} of the system, and $L$ and $S$ are operators determining the \emph{coupling} of the system to the field, with $S$ unitary.
%where the observable $\mathcal{H}$ is called the \emph{Hamiltonian} and ${\pmb i}$ denotes the imaginary unit. 
The evolution of $\psi$ is equivalent to the evolution of the observable $X$ given by \begdi X(t)=U^\ast (t) (X\otimes I)\, U(t), \enddi whose evolution is referred as %former is referred as the \emph{Schr\"odinger picture}, while the latter is 
the \emph{Heisenberg picture} while the one for $\psi$ is known as the \emph{Schr\"odinger picture}. This paper exclusively takes the point of view of the Heisenberg picture. The quantum stochastic calculus in \cite{Parthasarathy_92} allows to express the Heisenberg picture evolution of an scalar operator $X$ interacting with a boson field as
\begdi %eq \label{eq:general_evolution}
\begin{split}
dX  = &   \, (S^\ast X S-X)\,d \Lambda_w + {\cal L}(X) \,dt  + S^\ast[X,L] \,dW^\ast \\
& {} + [L^\ast,X]S \, dW,
\end{split}
\enddi %eq
where ${\cal L}(X)$ is the Lindblad operator defined as
\begdi %eq \label{eq:Linblad_operator}
{\cal L}(X)=-{\pmb i}[X,\mathcal{H}] +\frac{1}{2}\left(L^\ast[X,L] + [L^\ast,X]L\right). 
\enddi %eq
The output field is given by \begdi Y(t)=U(t)^{\ast}W(t)U(t), \enddi which amounts to %the QSDE %is 
\begdi
dY = L dt + S dW. 
\enddi
In summary, the dynamics of an open quantum system is uniquely determined by the triple $(S,L,\mathcal{H})$. Hereafter, the operator $S$ is assumed to be the identity operator ($S=I$). If on the other hand one consider $n_w$ interacting boson fields then the evolution equation is written as 
\begdi
dX= \mathcal{L}(X)\,dt + dW^\dagger\, [X,L]  + [L^\dagger,X]\, dW,
\enddi
where $[X,dW]=[X,dW^\dagger]^T=0$, $L=(L_1, \hdots , L_s )^T$,
\begdi
\mathcal{L}(X) \triangleq -\pmb{i} [X,H] +\frac{1}{2} \left(L^\dagger \, [X,L] + [L^\dagger,X]\,L \right),
\enddi
\begdi
dW=\left(\begin{array}{c}
    dW_1 \\ \vdots \\ dW_{n_w}
   \end{array}\right) \;\; {\rm and} \;\;
dW^\dagger = \left( dW^\ast_1, \cdots , dW^\ast_{n_w} \right).
\enddi

Consider the vector of operators $x = (x_1,\hdots, x_s)^T$. By stacking (column-wise) the scalar evolutions for each $x_i$, it follows that the Heisenberg evolution equation is

\begin{align} \label{eq:general_evolution_vector}
\nonumber \left(\begin{array}{c} dx_1  \\ \vdots \\ dx_s \end{array}\right) = & \; \left(\begin{array}{c} \mathcal{L}(x_1)  \\ \vdots \\ \mathcal{L}(x_s) \end{array}\right) \,dt +  \left(\begin{array}{c} [x_1,L^T]  \\ \vdots \\ {[x_s,L^T]} \end{array}\right)  \, dW^\# \\
\nonumber    & \, + \left(\begin{array}{c} [L^\dagger, x_1]  \\ \vdots \\ {[L^\dagger, x_s]} \end{array}\right)\, dW \\
\nonumber  = & \, \left(\begin{array}{c} \mathcal{L}(x_1)  \\ \vdots \\ \mathcal{L}(x_s) \end{array}\right) \,dt +  \left(\begin{array}{c} [x_1,L^T]  \\ \vdots \\ {[x_s,L^T]} \end{array}\right)  \, dW^\# \\
\nonumber    & \,  - \left(\begin{array}{c} [x_1 , L^\dagger]  \\ \vdots \\ {[x_s , L^\dagger]} \end{array}\right)\, dW \\
           = & \; \mathcal{L}(x) \,dt +  [x,L^T]  \, dW^\# - [x,L^\dagger]\, dW, 
\end{align}
where 
\begin{align}  \label{eq:Lindblad_operator_vector}
\nonumber \lefteqn{\hspace*{-0.2in} \mathcal{L}(x) } \\ 
\nonumber   = & \, -\pmb{i} [x,H] \\
\nonumber     & \, + \frac{1}{2} \left( \left(\begin{array}{c} L^\dagger \, [x_1,L]   \\ \vdots \\ L^\dagger \, {[x_s,L]} \end{array}\right) +  \left(\begin{array}{c} [L^\dagger,x_1]\,L  \\ \vdots \\ {[L^\dagger,x_s]}\,L \end{array}\right) \right) \\
\nonumber   = & \, -\pmb{i} [x,H] - \frac{1}{2} \left( \left(L^\dagger \! \left(  [L, x_1], \cdots , {[L,x_s]} \rule{0in}{0.15in}\right) \right)^T \right) \\
\nonumber     &\, - \frac{1}{2}[x, L^\dagger]\,L  \\
\nonumber   = & \, -\pmb{i} [x,H] + \frac{1}{2} \left( -\left(L^\dagger \,  [L, x^T]   \right)^T   -[x, L^\dagger]\,L  \right) \\
            = & \, -\pmb{i} [x,H] + \frac{1}{2} \left( \left(L^\dagger \,  [x, L^T]^T   \right)^T   +[ L^\#,x^T]^T\,L  \right).
\end{align}
% This amounts to the following Heisenberg evolution equation
% \begeq \label{eq:general_evolution_vector}
% dx= \mathcal{L}(x)\,dt +  [x,L^T]\,dW^\#  - [x,L^\dagger]\, dW,
% \endeq
% where 
% \begin{align} \label{eq:Lindblad_operator_vector}
% \nonumber \mathcal{L}(x) \triangleq & -\pmb{i} [x,H] \\
% & \,+\frac{1}{2} \left(\left( L^\dagger \, \left[x,L^T\right]^T \right)^T + \left[L^\#,x^T \right]^T \,L \right).
% \end{align}

It is customary to express QSDEs in terms of its interaction with quadrature fields. The quadrature fields are given by the transformation
\begeq \label{eq:quadrature_transformation}
\left(\begin{array}{c}
\bar{W}_1 \\ \bar{W}_2
\end{array}\right)= \left(\begin{array}{cc}
I_{n_w} & I_{n_w} \\ -\pmb{i}I_{n_w} & \pmb{i}I_{n_w} 
\end{array}\right) 
\left(\begin{array}{c}
W \\ W^\#
\end{array}\right),
\endeq
where the operators $\bar{W}_1$ and $\bar{W}_2$ are now self-adjoint, and $I_{n_w}$ denotes the identity matrix of dimension $n_w$. In \cite{Hudson_parthasarathy_84}, the It\^o table for $W$ and $W^\dagger$ is  
\begdi
\left(\begin{array}{c}
dW \\ dW^\#
\end{array}\right)    \left(\begin{array}{cc}
dW & dW^\#
\end{array}\right)= \left(\begin{array}{cc}
0 & I_{n_w} \\ 0 & 0 
\end{array}\right) dt,
\enddi
which in terms of the quadrature fields is
\begdi
\left(\begin{array}{c}
d\bar{W}_1 \\ d\bar{W}_2
\end{array}\right)    \left(\begin{array}{cc}
d\bar{W}_1 & d\bar{W}_2
\end{array}\right)= \left(\begin{array}{cc}
I_{n_w} & \pmb{i}I_{n_w} \\ -\pmb{i}I_{n_w} & I_{n_w} 
\end{array}\right) dt.
\enddi
Thus,
\begin{align*} %\label{eq:general_evolution_vector}
dx = & \, \mathcal{L}(x)\,dt \\
& \, +  \frac{1}{2}\left([x,L^T], -[x,L^\dagger]\right)\left(\begin{array}{cc}
I_{n_w} & -\pmb{i}I_{n_w} \\ I_{n_w} & \pmb{i}I_{n_w} 
\end{array}\right) \left(\begin{array}{c} d\bar{W}_1 \\ d\bar{W}_2 \end{array}\right) \\%\,dW^\#  - [x,L^\dagger]\, dW,
= & \, \mathcal{L}(x)\,dt + \frac{1}{2}\left([x,L^T]-[x,L^\dagger] \right)d\bar{W}_1  \\
& \,  -\frac{\pmb{i}}{2}\left([x,L^T]+[x,L^\dagger] \right)d\bar{W}_2
\end{align*}
%The output field is also transformed into its quadrature form as
The quadrature form of the output fields is obtained from the quadrature transformation
\vspace*{0.05in}\begdi
\left(\begin{array}{c}
\bar{Y}_1 \\ \bar{Y}_2
\end{array}\right)= \left(\begin{array}{cc}
I_{n_w} & I_{n_w} \\ -\pmb{i}I_{n_w} & \pmb{i}I_{n_w} 
\end{array}\right) 
\left(\begin{array}{c}
Y \\ Y^\dagger
\end{array}\right), 
\enddi
which gives 
\begdi % \label{eq:bilinear_system_output}
\left(\begin{array}{c}
d\bar{Y}_1 \\ d\bar{Y}_2
\end{array}\right)= \left(\begin{array}{c}
L+L^\# \\ \pmb{i}(L^\#-L) 
\end{array}\right) \,dt + \left(\begin{array}{c}
d\bar{W}_1 \\ d\bar{W}_2
\end{array}\right).
\enddi
%where $C_1=C + C^\#$ and $C_2=\pmb{i}(C^\# - C)$ are obviously real matrices.

The main focus of this paper is on the dynamics of open multi-level quantum systems interacting with $n_w$ Boson quantum fields. Such systems evolve with respect to the group $SU(n)$. The algebra of $SU(n)$ has been extensively studied since the $1950's$ to the point that it is an standard topic in quantum mechanics when studying multi-level systems \cite{Mahler-Weberrus_98,Macfarlane-Sudbery-Weisz_68}. To particularize the framework presented in the previous paragraph for system \rref{eq:general_evolution_vector} evolving on $SU(n)$, consider the Hilbert space $\mathfrak{H}=\C^n$ and let $\ket{j}$ with $j=1,\hdots,n$ be eigenvectors spanning $\mathfrak{H}$. A \emph{projection operator} ${P}_{kl}$ is defined as the outer product
\begdi
{P}_{k,l} = \ket{k} \bra{l},
\enddi
where $k,l=1,\hdots,n$. It is a well-known fact that any operator defined in $\mathfrak{H}$ can be obtained in term of these $n^2$ projection operators. Specifically the generators of $SU(n)$ are constructed as follows
  \begin{align*}
u_{jk} & = P_{j,k}+P_{k,j},\\  
v_{jk} & = \pmb{i}\left(P_{j,k}-P_{k,j} \right),\\
w_l    & = -\sqrt{\frac{2}{l(l+1)}} \left( \sum_{s=1}^k P_{s,s} - k P_{l+1,l+1} \right) 
\end{align*}
for $1\le j < k \le n$, $1 \le l \le n-1$. The $(n^2-n)/2$ symmetric matrices $u_{jk}$, the $(n^2-n)/2$ antisymmetric matrices $v_{jk}$ and the $n-1$ mutually commutative matrices $w_l$ together form the set $\{\lambda_1, \hdots, \lambda_{n^2-1} \}$ of generators of $SU(n)$. Hereafter define $s = n^2-1$. Their commutation and anticommutation relations are
\begin{align*}
[\lambda_i,\lambda_j] & = 2 \pmb{i} \sum_{k=1}^{s} f_{ijk} \lambda_{k}, \\
\{\lambda_i,\lambda_j \} & = \frac{4}{n} \delta_{ij} + 2 \sum_{k=1}^{s} d_{ijk} \lambda_{k}. 
\end{align*}
Thus, the product $\lambda_i \lambda_j$ can be easily computed as
\begin{align} \label{eq:product_Gell-mann_matrices}
\nonumber \lambda_i \lambda_j & = \frac{1}{2} \left( [\lambda_i,\lambda_j] + \{ \lambda_i,\lambda_j \} \right) \\ 
                    & = \frac{2}{n} \delta_{ij} + \sum_{k=1}^{s} \left( \pmb{i} f_{ijk} +d_{ijk}  \right)\lambda_k.
\end{align}
where the real completely antisymmetric tensor $f_{ijk}$ and the real completely symmetric tensor $d_{ijk}$ are called the \emph{structure constants} of $SU(n)$, and $\delta_{ij}$ is the Kronecker delta. The tensors $f_{ijk}$ and $d_{ijk}$ satisfy 
\begeq \label{eq:property_ff_tensor}
f_{ilm}f_{mjk}+f_{jlm}f_{imk}+f_{klm}f_{ijm}=0,
\endeq
\begeq \label{eq:property_fd_tensor}
f_{ilm}d_{mjk}+f_{jlm}d_{imk}+f_{klm}d_{ijm}=0,
\endeq
% \begeq \label{eq:property_fd_tensor}
% (\lambda_i)_{\alpha\beta} (\lambda_i)_{\gamma\rho} = 2 \delta_{\alpha\rho} \delta_{\beta\gamma} -\frac{2}{n}\delta_{\alpha\beta} \delta_{\gamma\rho}. 
% \endeq
% 
% \begdi
% \sum_{m=1}^s f_{ijm} f_{klm}  = \frac{2}{n} ( \delta_{ik} \delta_{jl} - \delta_{il} \delta_{jk} ) + (d_{ikm} d_{jlm}- d_{jkm} d_{ilm}). 
% \enddi
and
\begeq \label{eq:property_ff_delta}
\sum_{m,k=1}^s f_{im k} f_{jm k}  = n  \delta_{ij}. 
\endeq 

The procedure of how to construct the generalized Gell-Mann matrices shows that only the dimension of the group $SU(n)$ is necessary to express all the components of the group algebra, i.e., once $n$ is given then the generators and structure constants $f$ and $d$ are fixed. The vector of system variables for \rref{eq:general_evolution_vector} is 
\begdi 
x=\left(\begin{array}{cc}
         x_1 \\ \vdots \\ x_s
        \end{array} \right) \triangleq \left(\begin{array}{cc}
         \hat{\lambda}_1 \\ \vdots \\ \hat{\lambda}_s
        \end{array} \right),
\enddi where $\hat{\lambda}_1, \ldots, \hat{\lambda}_s$ are spin operators. Given that these operators are self-adjoint, the vector of operators $x$ satisfies $x=x^\#$.~In particular, a self-adjoint operator $\hat{\lambda}$ in $\mathfrak{H}$ is spanned by the generalized Gell-Mann matrices \cite{Mahler-Weberrus_98}, i.e.,
\begdi
\hat{\lambda}=\frac{1}{n} \alpha_0 + \frac{1}{2}\sum_{i=0}^s \alpha_i \lambda_i,
\enddi
where $\alpha_0=\Tr(\hat{\lambda})$, $\alpha_i=\Tr(\hat{\lambda}\lambda_i)$. Thus, $\alpha_0$ and $(\alpha_1,\ldots,\alpha_s)^T\in \C^3$ determine uniquely the operator $\hat{\lambda}$ with respect to a given basis in $\C^n$. The initial value of the system variables can be set to $x(0)=(\lambda_1,\ldots,\lambda_s)$.

It is important to emphasize that any arbitrary polynomial of the components of $x$ evolving in $SU(n)$ is described by a linear combination of the generators of $SU(n)$, which includes the identity \cite{Mahler-Weberrus_98}. This fact can be appreciated in relation \rref{eq:product_Gell-mann_matrices} because equation \rref{eq:general_evolution_vector} implies that such commutation relations are preserved by the evolution in time. Thus, any Hamiltonian and coupling operators of polynomial type are representable as linear functions of $x$. Therefore, assuming linearity captures a large class of Hamiltonian and coupling operators without much loss of generality, i.e., the assumed Hamiltonian is ${\mathcal{H}}=\alpha x$ with $\alpha\in \re^{s}$, and the multiplicative coupling operator is of the form $L=\Lambda x$ with $\Lambda \in \C^{n_w \times s	}$. The reason why a coupling operator is called multiplicative is that they make the interacting fields to appear in \rref{eq:general_evolution_vector} as multiplicative quantum noise.

% It is customary to express QSDEs in terms of its interaction with quadrature fields. The quadrature fields are given by the transformation
% \begeq \label{eq:quadrature_transformation}
% \left(\begin{array}{c}
% \bar{W}_1 \\ \bar{W}_2
% \end{array}\right)= \left(\begin{array}{cc}
% I_{n_w} & I_{n_w} \\ -\pmb{i}I_{n_w} & \pmb{i}I_{n_w} 
% \end{array}\right) 
% \left(\begin{array}{c}
% W \\ W^\#
% \end{array}\right),
% \endeq
% where the operators $\bar{W}_1$ and $\bar{W}_2$ are now self-adjoint, and $I_{n_w}$ denotes the identity matrix of dimension $n_w$. In \cite{Hudson_parthasarathy_84}, the It\^o table for $W$ and $W^\dagger$ is  
% \begdi
% \left(\begin{array}{c}
% dW \\ dW^\#
% \end{array}\right)    \left(\begin{array}{cc}
% dW & dW^\#
% \end{array}\right)= \left(\begin{array}{cc}
% 0 & I_{n_w} \\ 0 & 0 
% \end{array}\right) dt,
% \enddi
% which in terms of the quadrature fields is
% \begdi
% \left(\begin{array}{c}
% d\bar{W}_1 \\ d\bar{W}_2
% \end{array}\right)    \left(\begin{array}{cc}
% d\bar{W}_1 & d\bar{W}_2
% \end{array}\right)= \left(\begin{array}{cc}
% I_{n_w} & \pmb{i}I_{n_w} \\ -\pmb{i}I_{n_w} & I_{n_w} 
% \end{array}\right) dt.
% \enddi
In general, the evolution of $x$ in quadrature form falls into a class of bilinear QSDEs expressed as
% \begin{align*}% \label{eq:bilinear_system}
% \lefteqn{dx  %= & \, A_0\,dt+A x \, dt +\sum_{k=1}^{n_w} \bar{B}_{1k}x\,dW^\#_{k} + \bar{B}_{2k}x\, d{W}_{k}\\
%  =  \, A_0\,dt+A x \, dt} & \\
% & \,+ \left( \bar{B}_{11}x , \cdots , \bar{B}_{1n_w}x , \bar{B}_{21}x  , \cdots , \bar{B}_{2n_w} x \right)\left(\begin{array}{c} dW^\# \\ dW \end{array}\right),
% \end{align*}
% whose quadrature form is
\begin{align} \label{eq:bilinear_system}
\nonumber \lefteqn{ \hspace*{-0.1in} dx %& = A_0\,dt+A x \, dt +\left(\bar{B}_{11}x, \cdots, \bar{B}_{1n_w}x, \bar{B}_{21}x , \cdots, \bar{B}_{2n_w} x \right) \left(\begin{array}{cc} I_{n_w} & -\pmb{i}I_{n_w} \\ I_{n_w} & \pmb{i}I_{n_w} \end{array}\right) \left(\begin{array}{c} d\bar{W}_1\\ d\bar{W}_2 \end{array}\right)\\
 = A_0\,dt+A x \, dt} \\ %+\sum_{k=1}^{n_w} \left( {B}_{1k}x\,d\bar{W}_{1k} + {B}_{2k}x\, d\bar{W}_{2k} \right)\\
& \hspace*{-0.1in}\,+ \left( B_{11}x , \cdots , B_{1n_w}x , B_{21}x  , \cdots , B_{2n_w} x \right)\left(\begin{array}{c} d\bar{W}_1 \\ d\bar{W_2} \end{array}\right),
\end{align}
where $A_0\in \re^s$ and $A,{B}_{1k}\triangleq \bar{B}_{1k}+\bar{B}_{2k},{B}_{2k}\triangleq \pmb{i}(\bar{B}_{2k}-\bar{B}_{1k}) \in \re^{s\times s}$, $k=1,\hdots, n_w$. The fact that all matrices in \rref{eq:bilinear_system} are real is due to the quadrature transformation \rref{eq:quadrature_transformation}. The quadrature output fields are  
% \begdi 
% dY = C  x\,dt+ \frac{1}{2} \left(d\bar{W}_1 + \pmb{i} d\bar{W}_2\right),
% \enddi
whose quadrature form is %can be obtained from the transformation
% \vspace*{0.05in}\begdi
% \left(\begin{array}{c}
% \bar{Y}_1 \\ \bar{Y}_2
% \end{array}\right)= \left(\begin{array}{cc}
% I_{n_w} & I_{n_w} \\ -\pmb{i}I_{n_w} & \pmb{i}I_{n_w} 
% \end{array}\right) 
% \left(\begin{array}{c}
% Y \\ Y^\dagger
% \end{array}\right). 
% \enddi
% Thus, 
\begeq \label{eq:bilinear_system_output}
\left(\begin{array}{c}
d\bar{Y}_1 \\ d\bar{Y}_2
\end{array}\right)= \left(\begin{array}{c}
C_1 \\ C_2 
\end{array}\right) x \,dt + %\left(\begin{array}{cc}
% I_{n_w} & 0 \\ 0 & I_{n_w} 
% \end{array}\right) 
\left(\begin{array}{c}
d\bar{W}_1 \\ d\bar{W}_2
\end{array}\right),
\endeq
where $C\in \C^{n_w \times s}$, $C_1 \triangleq C + C^\# $ and $C_2 \triangleq \pmb{i}(C^\# - C)$. Note that the quadrature transformation makes $C_1$ and $C_2$ to be real matrices.

The objective of this paper is to determine conditions on the coefficients in \rref{eq:bilinear_system} and \rref{eq:bilinear_system_output} under which there exists the Hamiltonian and the coupling operator of the above form such that \rref{eq:bilinear_system}, \rref{eq:bilinear_system_output} can be written as \ref{eq:general_evolution_vector}; see Definition \ref{def:physical realizability} in Section \ref{sec:section4}.

\section{Notation and Preliminary Results} \label{sec:section3}

Define $F_i,D_i \in \re^{s\times s}$, $i\in \{1,\hdots, s\}$, such that their $(j,k)$ component is $({F}_i)_{jk} = f_{ijk}$ and $(D_i)_{jk} = d_{ijk} $, respectively. In particular, the set $\{-\pmb{i} {F}_1, \hdots, -\pmb{i} {F}_{s}  \}$ is the adjoint representation of $SU(n)$. In \cite{Macfarlane-Sudbery-Weisz_68,Kaplan-Resnikoff_67}, identities \rref{eq:property_ff_tensor} and \rref{eq:property_fd_tensor} were employed to obtain the following useful relationships
\begin{align}
%\label{eq:D_F_1} \frac{n}{2} D_i & = -\sum_{k}^{s} d_{ijk} F_j F_k \\
%\label{eq:D_F_2}\frac{n}{2} D_j & = -\sum_{k}^{s} D_k F_j F_k \\
\label{eq:D_F_3}      [F_i,F_j] & = - \sum_{k}^{s} f_{ijk} F_k \\
\label{eq:D_F_4}      [F_i,D_j] & = - \sum_{k}^{s} f_{ijk} D_k \\
\label{eq:D_F_5}      F_i D_j + F_jD_i & =  \sum_{k}^{s} d_{ijk} F_k \\
\label{eq:D_F_6}      D_i F_j + D_jF_i & =  \sum_{k}^{s} d_{ijk} F_k.
\end{align}
\begde
Let $\beta\in \C^s$. The linear mappings $\Theta^-,\Theta^+: \C^s \rightarrow \C^{s\times s}$ are defined as 
\begin{align}
\label{eq:Theta_definition} \Theta^-(\beta) &= \left(F_1^T \beta,\cdots , F_{s}^T \beta \right) = \left(\begin{array}{c}
         \beta^T F_1^T  \\ \vdots \\ \beta^T F_{s}^T  \end{array} \right),\\
\label{eq:Thetaplus_definition} \Theta^+(\beta) & = \left(D_1^T \beta,\cdots , D_{s}^T \beta \right) = \left(\begin{array}{c} \beta^T D_1^T  \\ \vdots \\ \beta^T D_{s}^T  \end{array} \right).
\end{align}  
\endde
\vspace*{0.05in}
In order to simplify the notation, if $\beta$ is an $s$-dimensional row then it will be understood hereafter that $\Theta^-(\beta) = \Theta^-(\beta^T)$ and $\Theta^+(\beta)=\Theta^+(\beta^T)$. In addition, these two matrix functions are used to express the commutation and anticommutation relations of the vector of operators $x$ in a compact form. That is,
\begin{align}
\label{eq:lambda_CCR}  [x,x^T]    & = 2 \pmb{i} \Theta^-(x),\\
\label{eq:lambda_antiCCR} \{x,x^T \} & = \frac{4}{n} I + 2 \Theta^+(x).
\end{align}
It is also important to notice that the nature of the $f$ and $d$-tensors make the matrices $\Theta^-(\beta)$ and $\Theta^+(\beta)$ be antisymmetric and symmetric, respectively. Consider now the \emph{stacking operator} $\vec: \C^{m\times n} \rightarrow \C^{m n}$ whose action on a matrix creates a column vector by stacking its columns below one another. With the help of $\vec$, the matrices $\Theta^-(\beta)$ and $\Theta^+(\beta)$ can be reorganized so that 
\begdi
\vec(\Theta^-(\beta))=\left(\begin{array}{c} \Theta^-_1(\beta) \\ \vdots \\ \Theta^-_s(\beta) \end{array} \right)= F \beta, 
\enddi
and
\begdi
\vec(\Theta^+(\beta))=\left(\begin{array}{c} \Theta^+_1(\beta) \\ \vdots \\ \Theta^+_s(\beta) \end{array} \right)= D \beta, 
\enddi
where $\beta\in \C^s$, $\Theta^-_i(\beta)=F^T_i\beta$, $F=\left(F_1 ,\cdots ,F_s  \right)^T$, $\Theta^+_i(\beta)=D_i\beta$ and $D=\left(D_1 ,\cdots ,D_s  \right)^T$. From \rref{eq:property_ff_delta}, $F$ satisfies
\begin{align*} \label{eq:E^T_by_E}
(F^T  F )_{i j} & = \left( \sum_{k=1}^s F_k F^T_k \right)_{ij}   \\
& =   -\sum_{k,m=1}^s (F_k)_{im} (F_k)_{mj}  \\
& =   -\sum_{k,m=1}^s f_{kim} f_{kmj}  \\
& =   \sum_{k,m=1}^s f_{imk} f_{jmk}  \\
& =  n \delta_{ij},
\end{align*}
which implies 
\begeq \label{eq:E^T_by_E}
F^T  F= n I.
\endeq

The properties of $\Theta^-$ and $\Theta^+$ are summarized in the next lemma. 
\begle \label{le:prop_thetas} Let $\beta,\gamma\in \C^s$. The mappings $\Theta^-$ and $\Theta^+$ satisfy
%\vspace*{-0.1in}
\renewcommand*\theenumi{{\roman{enumi}}}
\renewcommand*\labelenumi{$\theenumi.$}
\begin{enumerate}
\item \label{eq:Theta_1} $\displaystyle \Theta^-(\beta) \gamma = -\Theta^-(\gamma) \beta,$
\vspace*{0.05in}
\item \label{eq:Theta_2} $\displaystyle \Theta^+(\beta) \gamma = \Theta^+(\gamma) \beta,$
\vspace*{0.05in}
\item \label{eq:Theta_beta_beta} $\displaystyle \Theta^-(\beta)\beta = 0,$
\vspace*{0.05in}
\item \label{eq:Theta__minus_theta_minus} $\displaystyle \Theta^-\left(\Theta^-(\beta) \gamma \right) = [\Theta^-(\beta) ,\Theta^-(\gamma) ],$% - \Theta^-(\gamma) \Theta^-(\beta),$
\vspace*{0.05in}
\item \label{eq:Theta__minus_theta_plus} $\displaystyle \Theta^-\left(\Theta^+(\beta) \gamma \right) = \{\Theta^-(\beta) , \Theta^+(\gamma)\},$% + \Theta^-(\gamma) \Theta^+(\beta),$ 
\vspace*{0.05in}
  \item \label{eq:Theta__plus_theta_minus1} $\displaystyle \Theta^+\left(\Theta^-(\beta) \gamma \right) = [\Theta^+(\beta) , \Theta^-(\gamma)] = [\Theta^-(\beta) , \Theta^+(\gamma) ],$\vspace*{0.05in}%- \Theta^-(\gamma) \Theta^+(\beta),$
  % \vspace*{0.05in}
  % \item \label{eq:Theta__plus_theta_minus2} $\displaystyle \Theta^+\left(\Theta^-(\beta) \gamma \right) = \Theta^-(\beta) \Theta^+(\gamma) - \Theta^+(\gamma) \Theta^-(\beta).$\vspace*{0.05in}
  \end{enumerate}
 where the commutator of matrices is defined as usual, i.e., $[A,B]=AB-BA$ for $A,B\in \re^{s\times s}$.
 \endle\vspace*{0.05in}

  \begpr Using \rref{eq:Theta_definition}, one can decompose the left-hand-side of $\rref{eq:Theta_1}$ in terms of the matrices $F_i$ as  
  \begdi
  \Theta^-(\beta) \gamma = \left( \begin{array}{c}  \beta^T F_1^T \gamma    \\  \vdots  \\ \beta^T F_{s}^T \gamma  \end{array}
  \right)
  \enddi
  Every component is the written as
  \begdi
  \beta^T F_i^T \gamma %= - \beta^T\left( \begin{array}{c} \displaystyle \sum_{k=1}^{s} f_{i1k} \gamma_k  \\  \vdots  \\ \displaystyle  \sum_{k=1}^{s} f_{isk} \gamma_k  \end{array}  \right) 
  = - \left( \begin{array}{c}  \displaystyle \sum_{k,l=1}^{s} \beta_l f_{ilk} \gamma_k  \\  \vdots  \\  \displaystyle  \sum_{k,l=1}^{s} \beta_l f_{ilk} \gamma_k  \end{array} \right)
= \left( \begin{array}{c} \displaystyle   \sum_{k,l=1}^{s} \gamma_k f_{ikl} \beta_l  \\  \vdots  \\  \displaystyle  \sum_{k,l=1}^{s} \gamma_k f_{ikl} \beta_l   \end{array} \right),
\enddi
which implies $\beta^T F_i^T \gamma = -\gamma^T F_i^T \beta$. Therefore $\Theta^-(\beta) \gamma = -\Theta^-(\gamma) \beta$. Using \rref{eq:Thetaplus_definition}, a similar procedure is applied for identity $\rref{eq:Theta_2}$ in terms of the matrices $D_i$. The $i$-th component is given by 
\begdi
\beta^T D_i \gamma %=  \beta^T\left( \begin{array}{c}  \displaystyle   \sum_{k=1}^{s} d_{i1k} \gamma_k  \\  \vdots  \\  \displaystyle  \sum_{k=1}^{s} d_{isk} \gamma_k  \end{array}  \right) 
=  \left( \begin{array}{c}  \displaystyle  \sum_{k,l=1}^{s} \beta_l d_{ilk} \gamma_k  \\  \vdots  \\  \displaystyle  \sum_{k,l=1}^{s} \beta_l d_{ilk} \gamma_k  \end{array} \right)
=  \left( \begin{array}{c}  \displaystyle  \sum_{k,l=1}^{s} \gamma_k d_{ikl} \beta_l   \\  \vdots  \\  \displaystyle  \sum_{k,l=1}^{s} \gamma_k d_{ikl} \beta_l   \end{array} \right) ,
\enddi
which gives $\beta^T D_i \gamma = \gamma^T D_i \beta $. Thus, $\Theta^+(\beta) \gamma = \Theta^+(\gamma) \beta$. Identity $\rref{eq:Theta_beta_beta}$ is true since $f_{ijj} =0$ for all $i$ and $j$, and 
\begin{align*}
\sum_{k,l=1}^{s} \beta_l f_{ilk} \beta_k  = & \sum_{\substack{k,l=1 \\ k\neq l}}^{s} \beta_l f_{ilk} \beta_k \\
 = & \sum_{k < l} \beta_l f_{ilk} \beta_k  + \sum_{k > l} \beta_l f_{ilk} \beta_k \\
 = & \sum_{k < l} \beta_l f_{ilk} \beta_k  -\sum_{k < l} \beta_l f_{ilk} \beta_k\\
 = & 0,
\end{align*}
where the negative sign in the last summand was obtained because of the antisymmetry of $f_{ilk}$. The left-hand-side of identity $\rref{eq:Theta__minus_theta_minus}$ is decomposed as 
\begin{align*}
\lefteqn{\Theta^-\left(\Theta^-(\beta) \gamma \right)} \\
 & = \left(  F_1 \left(\begin{array}{c} \beta^T F_1 \gamma \\ \vdots \\  \beta^T F_s \gamma \end{array}\right), \cdots, F_s \left(\begin{array}{c} \beta^T F_1 \gamma \\ \vdots \\  \beta^T F_s \gamma \end{array}\right) \right) \\
& = \left( \begin{array}{ccc} \displaystyle \sum_{k=1}^{s} f_{11k} \beta^T F_k \gamma & \cdots & \displaystyle \sum_{k=1}^{s} f_{s1k} \beta^T F_k \gamma  \\  \vdots &  \ddots &  \vdots  \\ \displaystyle \sum_{k=1}^{s} f_{1sk} \beta^T F_k \gamma & \cdots & \displaystyle  \sum_{k=1}^{s} f_{ssk} \beta^T F_k \gamma   \end{array} \right).
\end{align*}
By \rref{eq:D_F_3}, the $(i,j)$ component of this matrix is 
\begin{align*}
\left( \Theta^-\left(\Theta^-(\beta) \gamma \right) \right)_{ij} & = - \sum_{k=1}^{s} f_{ijk} \beta^T F_k \gamma \\
& =  \beta^T \left(-\sum_{k=1}^{s} f_{ijk} F_k \right) \gamma \\
& =  \beta^T [F_i, F_j] \gamma \\
& =  \beta^T F_iF_j \gamma - \beta^T F_jF_i \gamma \\
& =  \beta^T F_i^TF_j^T \gamma - \beta^T F_j^TF_i^T \gamma \\
& =  \beta^T F_i^TF_j^T \gamma - \gamma^T F_i^TF_j^T \beta \\
& = \left( \Theta^-(\beta) \Theta^-(\gamma) - \Theta^-(\gamma) \Theta^-(\beta) \right)_{ij}\\
& = \left([ \Theta^-(\beta), \Theta^-(\gamma) ] \right)_{ij}.   
\end{align*}
Similarly, decomposing the $(i,j)$ component of the left-hand-side of $\rref{eq:Theta__minus_theta_plus}$ and using \rref{eq:D_F_4} gives
% \begin{align*}
% \lefteqn{ \Theta^-\left(\Theta^+(\beta) \gamma \right)} \\
% & = - \left( F_1 \left(\begin{array}{c} \beta^T D_1 \gamma \\ \vdots \\  \beta^T D_s \gamma \end{array}\right), \cdots, F_s \left(\begin{array}{c} \beta^T D_1 \gamma \\ \vdots \\  \beta^T D_s \gamma \end{array}\right) \right) \\
% & = - \left( \begin{array}{ccc} \displaystyle \sum_{k=1}^{s} f_{11k} \beta^T D_k \gamma & \cdots & \displaystyle \sum_{k=1}^{s} f_{s1k} \beta^T D_k \gamma  \\  \vdots &  \ddots &  \vdots  \\ \displaystyle \sum_{k=1}^{s} f_{1sk} \beta^T D_k \gamma & \cdots & \displaystyle  \sum_{k=1}^{s} f_{ssk} \beta^T D_k \gamma   \end{array} \right).
% \end{align*}
% By \rref{eq:D_F_4}, the $(i,j)$ component of this matrix is 
\begin{align*}
\left( \Theta^-\left(\Theta^+(\beta) \gamma \right) \right)_{ij} & =  \sum_{k=1}^{s} f_{ijk} \beta^T D_k \gamma \\
& = - \beta^T \left(-\sum_{k=1}^{s} f_{ijk} D_k \right) \gamma \\
& = - \beta^T [F_i, D_j] \gamma \\
& = - \beta^T F_i D_j \gamma + \beta^T D_j F_i \gamma \\
& =  \beta^T F_i^TD_j \gamma - \beta^T D_j F_i^T \gamma \\
& =  \beta^T F_i^T D_j \gamma + \gamma^T F_i^T D_j \beta \\
& = \left( \Theta^-(\beta) \Theta^+(\gamma) + \Theta^-(\gamma) \Theta^+(\beta) \right)_{ij}\\
& = \left( \{\Theta^-(\beta), \Theta^+(\gamma) \} \right)_{ij}.
\end{align*}
Again, decomposing the $(i,j)$ component of the left-hand-side of $\rref{eq:Theta__plus_theta_minus1}$ and using \rref{eq:D_F_5} gives  %in terms of $F_i$'s and $D_i$'s is
% \begin{align*}
% \lefteqn{ \Theta^-\left(\Theta^+(\beta) \gamma \right)} \\
% & = - \left( F_1 \left(\begin{array}{c} \beta^T D_1 \gamma \\ \vdots \\  \beta^T D_s \gamma \end{array}\right), \cdots, F_s \left(\begin{array}{c} \beta^T D_1 \gamma \\ \vdots \\  \beta^T D_s \gamma \end{array}\right) \right) \\
% & = - \left( \begin{array}{ccc} \displaystyle \sum_{k=1}^{s} f_{11k} \beta^T D_k \gamma & \cdots & \displaystyle \sum_{k=1}^{s} f_{s1k} \beta^T D_k \gamma  \\  \vdots &  \ddots &  \vdots  \\ \displaystyle \sum_{k=1}^{s} f_{1sk} \beta^T D_k \gamma & \cdots & \displaystyle  \sum_{k=1}^{s} f_{ssk} \beta^T D_k \gamma   \end{array} \right).
% \end{align*}
% By \rref{eq:D_F_4}, the $(i,j)$ component of this matrix is 
\begin{align*}
\left( \Theta^+\left(\Theta^-(\beta) \gamma \right) \right)_{ij} & =  -\sum_{k=1}^{s} d_{ijk} \beta^T F_k \gamma \\
& = - \beta^T \left(\sum_{k=1}^{s} d_{ijk} F_k \right) \gamma \\
& = - \beta^T (F_i D_j+F_jD_i ) \gamma \\
%& = - \beta^T F_i D_j \gamma + \beta^T D_j F_i \gamma \\
%& =  \beta^T D_i F_j^T \gamma - \beta^T D_jF_i \gamma \\
& =  \beta^T F_i^T D_j \gamma - \gamma^T D_i F^T_j \beta \\
& = \left( \Theta^-(\beta) \Theta^+(\gamma) - \Theta^+(\gamma) \Theta^-(\beta) \right)_{ij}\\
& = \left( [ \Theta^-(\beta) , \Theta^+(\gamma) ] \right)_{ij}.
\end{align*}
Finally, applying the same procedure but %to $\rref{eq:Theta__plus_theta_minus2}$ 
using \rref{eq:D_F_6} instead gives
% \begin{align*}
% \lefteqn{ \Theta^-\left(\Theta^+(\beta) \gamma \right)} \\
% & = - \left( F_1 \left(\begin{array}{c} \beta^T D_1 \gamma \\ \vdots \\  \beta^T D_s \gamma \end{array}\right), \cdots, F_s \left(\begin{array}{c} \beta^T D_1 \gamma \\ \vdots \\  \beta^T D_s \gamma \end{array}\right) \right) \\
% & = - \left( \begin{array}{ccc} \displaystyle \sum_{k=1}^{s} f_{11k} \beta^T D_k \gamma & \cdots & \displaystyle \sum_{k=1}^{s} f_{s1k} \beta^T D_k \gamma  \\  \vdots &  \ddots &  \vdots  \\ \displaystyle \sum_{k=1}^{s} f_{1sk} \beta^T D_k \gamma & \cdots & \displaystyle  \sum_{k=1}^{s} f_{ssk} \beta^T D_k \gamma   \end{array} \right).
% \end{align*}
% By \rref{eq:D_F_4}, the $(i,j)$ component of this matrix is 
\begin{align*}
\left( \Theta^+\left(\Theta^-(\beta) \gamma \right) \right)_{ij} & =  -\sum_{k=1}^{s} d_{ijk} \beta^T F_k \gamma \\
& = - \beta^T \left(\sum_{k=1}^{s} d_{ijk} F_k \right) \gamma \\
& = - \beta^T (D_i F_j+D_jF_i ) \gamma \\
%& = - \beta^T F_i D_j \gamma + \beta^T D_j F_i \gamma \\
%& =  \beta^T D_i F_j^T \gamma - \beta^T D_jF_i \gamma \\
& =  \beta^T D_i F_j^T \gamma - \gamma^T F_i^T D_j \beta \\
& = \left( \Theta^+(\beta) \Theta^-(\gamma) - \Theta^-(\gamma) \Theta^+(\beta) \right)_{ij}\\
& = \left( [ \Theta^+(\beta) , \Theta^-(\gamma) ] \right)_{ij},
\end{align*}
which completes the proof. 
\endpr

\begle \label{le:commutators_evolution_equation}
Let $A,B\in \C^{n_w\times s}$, and denote by $A_i,B_i$ their respective rows, $i=1,\hdots,n_w$. Then
\begin{subequations}\label{eqn:commutators_general}
    \begin{alignat}{2}
&[x,(Ax)^T]  =  \, -2\pmb{i} \left( \Theta^-(A_1)x ,  \cdots , \Theta^-(A_{n_w}) x   \right), \label{eq:commutator_1}\\
\nonumber & [x,(Ax)^T]\, Bx   =  \, -2 \pmb{i} \sum_{k=1}^{n_w} \left( \frac{2}{n} \Theta^-(A_k)  B_k^T  \right. \\
&\, \left. \hspace*{0.1in} \rule{0in}{0.2in} + \Theta^-(A_k) \Theta^+(B_k) x - \pmb{i} \Theta^-(A_k) \Theta^-( B_k) x   \right), \label{eq:commutator_2}\\
\nonumber & (Bx)^T[Ax,x^T] =  \, 2\pmb{i} \sum_{k=1}^{n_w} \left(\frac{2}{n}\Theta^-(A_k)B_k^T  \right. \\
&\, \left. \hspace*{0.1in} \rule{0in}{0.2in} + \Theta^-(A_k) \Theta^+(B_k) x + \pmb{i}\Theta^-(A_k)\Theta^-( B_k) x \right)^T\!\!\!. \label{eq:commutator_3}
      \end{alignat}
\end{subequations}
%\vspace*{0.05in}
\endle
\begpr
The goal is to rewrite \rref{eq:commutator_1} in terms of $[x,x^T]$ and $\{x,x^T\}$ in order to apply \rref{eq:lambda_CCR} and \rref{eq:lambda_antiCCR}. Then 
\begin{align*}
\lefteqn{\hspace*{-0.15in}[x,(Ax)^T]} \\
 = & \, \left(\begin{array}{ccc} x_1x^TA_1^T &  \cdots & x_1x^TA_{n_w}^T \\ \vdots & \ddots & \vdots \\ x_sx^TA_1^T &  \cdots & x_sx^TA_{n_w}^T \end{array}\right) \\
& \, - \left(\begin{array}{ccc} A_1 x x_1 &  \cdots & A_1 x x_1 \\ \vdots & \ddots & \vdots \\ A_{n_w}xx_1 &  \cdots & A_{n_w} x x_s \end{array}\right)^T \\
%      & = \left(\begin{array}{ccc} x_1x^TA_1^T &  \cdots & x_1x^TA_{n_w}^T \\ \vdots & \ddots & \vdots \\ x_sx^TA_1^T &  \cdots & x_sx^TA_{n_w}^T \end{array}\right) - \left(\begin{array}{ccc} A_1 x x_1 &  \cdots &  A_{n_w}xx_1 \\ \vdots & \ddots & \vdots \\ A_1 x x_1 &  \cdots & A_{n_w} x x_s \end{array}\right) \\
%      & = \left( \left(\begin{array}{c} x_1A_1 x - ((A_1x) x_1 )^T \\  \vdots \\ x_s A_{1} x - ((A_{1}x) x_s )^T  \end{array}\right) , \cdots , \left(\begin{array}{c} x_1A_{n_w} x - ((A_{n_w}x) x_1 )^T \\  \vdots \\ x_s A_{n_w} x - ((A_{n_w}x) x_s )^T  \end{array}\right) \right) \\
       = & \,\left( x x^TA_1^T - ((A_1x) x^T )^T ,  \cdots , x x^TA_{n_w}^T - ((A_1x) x^T )^T \right) \\
       = & \,\left( x x^TA_1^T - (x x^T)^T A_1^T ,  \cdots , x x^TA_{n_w}^T - (x x^T)^T A_{n_w}^T \right) \\
       = & \,\left( [x ,x^T] A_1^T ,  \cdots , [x, x^T] A_{n_w}^T  \right) \\
       = & \, 2\pmb{i} \left( \Theta^-(x)A_1^T ,  \cdots , \Theta^-(x) A_{n_w}^T  \right).
\end{align*}
Thus, Lemma \ref{le:prop_thetas} gives  
\begdi
[x,(Ax)^T] = - 2\pmb{i} \left( \Theta^-(A_1)x ,  \cdots , \Theta^-( A_{n_w}) x  \right).
\enddi
For \rref{eq:commutator_2}, note that the scalar operator $B_i x$ commutes with $\Theta^-(A_j)$ for any $i$ and $j$. Recall that \begdi xx^T=\frac{1}{2}([x,x^T]+\{x,x^T\}).\enddi It then follows that\vspace*{0.05in}
\begin{align*}
\lefteqn{[x,(Ax)^T]\, Bx} \\
= & \, - 2\pmb{i} \left( \Theta^-(A_1)x ,  \cdots , \Theta^-( A_{n_w}) x  \right) \left(\begin{array}{c} B_1 x \\  \vdots \\ B_{n_w} x  \end{array}\right)\\
= & \, - 2\pmb{i}  \left(\begin{array}{ccc} A_1 F_1^T x &  \cdots & A_{n_w} F_1^T x \\ \vdots & \ddots & \vdots \\ A_1 F_s^T x &  \cdots & A_{n_w} F_s^T x \end{array}\right) \left(\begin{array}{c} B_1 x \\  \vdots \\ B_{n_w} x  \end{array}\right)\\
%= & \, -2 \pmb{i} \sum_{k=1}^{n_w} \left(\begin{array}{c} A_k F_1^T x B_k x \\  \vdots \\ A_k F_s^T x  B_{k} x  \end{array}\right) \\
= & \, -2 \pmb{i} \sum_{k=1}^{n_w} \left(\begin{array}{c} A_k F_1^T xx^T B_k^T \\  \vdots \\ A_k F_s^T x x^T B_k^T   \end{array}\right) \\
= & \, -2 \pmb{i} \sum_{k=1}^{n_w} \left(\begin{array}{c} A_k F_1^T \left(\frac{2}{n}I_s +\Theta^+(x) +\pmb{i}\Theta^-(x) \right) B_k^T \\  \vdots \\ A_k F_s^T \left(\frac{2}{n}I_s +\Theta^+(x) +\pmb{i}\Theta^-(x) \right) B_k^T   \end{array}\right) \\
%= & \, -2 \pmb{i} \sum_{k=1}^{n_w} \left( \frac{2}{n} \left(\begin{array}{c} A_k F_1^T  B_k^T \\  \vdots \\ A_k F_s^T  B_k^T   \end{array}\right) + \left(\begin{array}{c} A_k F_1^T \Theta^+(x)  B_k^T \\  \vdots \\ A_k F_s^T \Theta^+(x) B_k^T   \end{array}\right) + \pmb{i} \left(\begin{array}{c} A_k F_1^T \Theta^-(x)  B_k^T \\  \vdots \\ A_k F_s^T \Theta^-(x)  B_k^T   \end{array}\right)  \right) \\
%= & \, -2 \pmb{i} \sum_{k=1}^{n_w} \left( \frac{2}{n} \left(\begin{array}{c} \Theta_1^-(A_k)^T  B_k^T \\  \vdots \\ \Theta_s^-(A_k)^T  B_k^T   \end{array}\right)  \right.\\
% & \, \left. + \left(\begin{array}{c} \Theta_1^-(A_k)^T \Theta^+(x)  B_k^T \\  \vdots \\ \Theta_s^-(A_k)^T \Theta^+(x) B_k^T   \end{array}\right) + \pmb{i} \left(\begin{array}{c} \Theta_1^-(A_k)^T \Theta^-(x)  B_k^T \\  \vdots \\ \Theta_s^-(A_k)^T \Theta^-(x)  B_k^T   \end{array}\right)  \right) \\
= & \, -2 \pmb{i} \sum_{k=1}^{n_w} \left( \frac{2}{n} \Theta^-(A_k)  B_k^T  + \Theta^-(A_k) \Theta^+(x)  B_k^T \right. \\
 & \, \left. \rule{0in}{0.2in} \hspace*{0.6in} + \pmb{i} \Theta^-(A_k) \Theta^-(x)  B_k^T   \right) \\
= & \, -2 \pmb{i} \sum_{k=1}^{n_w} \left( \frac{2}{n} \Theta^-(A_k)  B_k^T  + \Theta^-(A_k) \Theta^+(B_k) x \right. \\
 & \, \left. \rule{0in}{0.2in} \hspace*{0.6in} - \pmb{i} \Theta^-(A_k) \Theta^-( B_k) x   \right).
\end{align*}
Similarly for \rref{eq:commutator_3}, one has, applying \rref{eq:commutator_vector_transpose}, that %since $[x,(Ax)^T]^T=-[Ax,x^T]$, one has that
\vspace*{0.05in}\begin{align*}
\lefteqn{(Bx)^T[Ax,x^T]} \\ %& =  2\pmb{i} \left( B_1 x ,  \cdots , B_{n_w} x  \right)  \left(\begin{array}{ccc} A_1 F_1^T x &  \cdots & A_{1} F_s^T x \\ \vdots & \ddots & \vdots \\ A_{n_w} F_1^T x &  \cdots & A_{n_w} F_s^T x \end{array}\right) \\
       & = 2\pmb{i} \sum_{k=1}^{n_w} \left( B_k x A_k F_1^T x ,  \cdots , B_k x A_k F_s^T x  \right) \\
       & = 2\pmb{i} \sum_{k=1}^{n_w} \left( B_k x \, \Theta_1^-(A_k) x ,  \cdots , B_k x \, \Theta_s^-(A_k) x  \right) \\
%       & = 2\pmb{i} \sum_{k=1}^{n_w} \left( \Theta_1^-(A_k) (B_k x) x ,  \cdots , \Theta_s^-(A_k) (B_k x) x  \right) \\
       & = 2\pmb{i} \sum_{k=1}^{n_w} \left( \Theta_1^-(A_k)  (x x^T)^T B_k^T , \cdots , \Theta_s^-(A_k) (x x^T)^T B_k^T \right) \\
%       & = 2\pmb{i} \sum_{k=1}^{n_w} \left( \Theta_1^-(A_k)  \left(\frac{2}{n}I_s +\Theta^+(x) - \pmb{i}\Theta^-(x) \right) B_k^T ,  \cdots , \Theta_s^-(A_k) \left(\frac{2}{n}I_s +\Theta^+(x) -\pmb{i}\Theta^-(x) \right) B_k^T \right) \\
       & = 2\pmb{i} \sum_{k=1}^{n_w} \left(\Theta^-(A_k)  \left(\frac{2}{n}I_s +\Theta^+(x) - \pmb{i}\Theta^-(x) \right) B_k^T \right)^T \\
       & = 2\pmb{i} \sum_{k=1}^{n_w} \left(\frac{2}{n}\Theta^-(A_k)B_k^T  + \Theta^-(A_k) \Theta^+(x)B_k^T \right. \\
 & \, \left. \rule{0in}{0.2in} \hspace*{0.6in} - \pmb{i}\Theta^-(A_k)\Theta^-(x) B_k^T \right)^T \\
       & = 2\pmb{i} \sum_{k=1}^{n_w} \left(\frac{2}{n}\Theta^-(A_k)B_k^T  + \Theta^-(A_k) \Theta^+(B_k) x \right. \\
 & \, \left. \rule{0in}{0.2in} \hspace*{0.6in} + \pmb{i}\Theta^-(A_k)\Theta^-( B_k) x \right)^T\!\!\!.\vspace*{-0.1in}
\end{align*} \endpr

The explicit computation of the vector fields in \rref{eq:general_evolution_vector} is given in the next lemma. 

\vspace*{0.1in}\begle \label{le:commutators_evaluated}  The component coefficients of equations \rref{eq:general_evolution_vector} and \rref{eq:Lindblad_operator_vector} are 
\begin{subequations}\label{eqn:commutators_for_evolution}
    \begin{alignat}{2}
\hspace*{-0.15in}  &[x,\mathcal{H}]  =  \,-2{\pmb i}\Theta^-(\alpha)x, \label{eq:commutator_xH}\\
\hspace*{-0.15in}  &[x,L^T]  =  \, -2{\pmb i} \left( \Theta^-(\Lambda_1)x ,\cdots, \Theta^-(\Lambda_{n_w})x \right), \label{eq:commutator_xL}\\
\hspace*{-0.15in}  &[x,L^\dagger] =  \, -2{\pmb i}\left( \Theta^-(\Lambda_1^{\#})x, \cdots , \Theta^-(\Lambda_{n_w}^{\#})x \right),  \label{eq:commutator_xLT}\\
\nonumber \hspace*{-0.15in}  &[L^\#,x^T]\, L  =  \, \sum_{k=1}^{n_w} \left( \frac{4 \pmb{i}}{n} \Theta^-(\Lambda^\#_k)  \Lambda_k^T \right. \\
& \, \left. \rule{0in}{0.2in}  + 2 \pmb{i}\Theta^-(\Lambda^\#_k) \Theta^+(\Lambda_k) x + 2\Theta^-(\Lambda^\#_k) \Theta^-( \Lambda_k) x   \right),\label{eq:commutator_LTxL}\\
\nonumber \hspace*{-0.15in}  &\left(L^\dagger \, \left[x,L^T\right]^T \right)^T \!\!\! =  \,  \sum_{k=1}^{n_w} \left( \frac{4 \pmb{i}}{n} \Theta^-(\Lambda^\#_k)  \Lambda^T_k \right. \\
 & \, \left. \rule{0in}{0.2in}  - 2 \pmb{i} \Theta^-(\Lambda_k) \Theta^+(\Lambda^\#_k) x + 2 \Theta^-(\Lambda_k) \Theta^-( \Lambda^\#_k) x   \right). \label{eq:commutator_xLTL}
    \end{alignat}
\end{subequations}
\endle
\begpr  Commutators \rref{eq:commutator_xH}-\rref{eq:commutator_xLT} follow directly from \rref{eq:commutator_1}. Commutator \rref{eq:commutator_LTxL} is computed out of \rref{eq:commutator_2} as
\begin{align*}
[L^\#,x^T]\, L & = - [x,L^\dagger] \, L \\
& =   \sum_{k=1}^{n_w} \left( \frac{4 \pmb{i}}{n} \Theta^-(\Lambda^\#_k)  \Lambda_k^T  + 2 \pmb{i}\Theta^-(\Lambda^\#_k) \Theta^+(\Lambda_k) x \right. \\
 & \, \left. \rule{0in}{0.2in} \hspace*{0.45in}+ 2\Theta^-(\Lambda^\#_k) \Theta^-( \Lambda_k) x   \right).
\end{align*}
Finally, commutator \rref{eq:commutator_xLTL} is obtained using \rref{eq:commutator_3} as
\begin{align*}
\left(L^\dagger \, \left[x,L^T\right]^T \right)^T = & \, - \left( L^\dagger \, \left[L^\#,x^T \right]  \right)^T  \\
 = & \,  \sum_{k=1}^{n_w} \left( \frac{4 \pmb{i}}{n} \Theta^-(\Lambda^\#_k)  \Lambda^T_k   \right. \\
 & \, \left. \rule{0in}{0.2in} \hspace*{-0.6in} - 2 \pmb{i} \Theta^-(\Lambda_k) \Theta^+(\Lambda^\#_k) x + 2 \Theta^-(\Lambda_k) \Theta^-( \Lambda^\#_k) x   \right)
\end{align*}
\endpr

\vspace*{-0.1in}

\section{Physical Realizability}   \label{sec:section4}

% In an environment where the classical laws of physics apply, standard control techniques such as optimization or a Lyapunov procedures do not worry in general of the nature of the controller they synthesized. In other words, their implementation is always possible since the physics behind them still hold. However, if one desires to implement a controller that obeys the laws imposed by quantum mechanics (quantum coherent control), then such a task is not so easily achieved unless an explicit characterization of that laws is given in terms of the control system vector fields. This is exactly the purpose for introducing the concept of a \emph{physically realizable} system in the next definition.  

The main contribution of the paper are given in this section. First, physical realizability is introduced next %definition%concept of physical

\begde \label{def:physical realizability} A system described by equations \rref{eq:bilinear_system} and \rref{eq:bilinear_system_output} is said to be {physically realizable} if there exist $\mathcal{H}$ and $L$ such that \rref{eq:bilinear_system} can be written as in \rref{eq:general_evolution_vector}.
\endde

The explicit form of matrices $A_0,A, B_{1k}, B_{2k},C_1$ and $C_2$ in terms of the Hamiltonian and coupling operator is given next.
\begth \label{th:physical_realizability_def} Let $\mathcal{H}=\alpha x$, with
$\alpha^T \in \re^s$, and $L=\Lambda x$, with $\Lambda \in \C^{n_w \times s}$. Then
\begin{subequations}\label{eqn:Spin_system_matrices}
    \begin{align}
      A_0 &= \frac{4 \pmb{i}}{n} \sum_{k=1}^{n_w}\Theta^-(\Lambda^\#_k)\Lambda_k^T, \label{subeqn:F0}\\
      A   &= - 2 \Theta^-(\alpha) + \sum_{k=1}^{n_w} \left(R_k -\pmb{i} Q_k \right), \label{subeqn:F}\\
      {B}_{1k} &= \Theta^-\left({\pmb i}(\Lambda^\#_k - \Lambda_k)\right), \label{subeqn:G1}\\ 
      {B}_{2k} &= \Theta^-(\Lambda_k + \Lambda^\#_k), \label{subeqn:G2} \\
      C_1   &=  \Lambda+\Lambda^\#, \label{subeqn:H12} \\
      C_2   &=  \pmb{i} \left( \Lambda^\#-\Lambda \right), \label{subeqn:H2}
   \end{align}
  \end{subequations}
where \begdi R_k\triangleq \Theta^-(\Lambda_k)\Theta^-(\Lambda^\#_k) + \Theta^-(\Lambda^\#_k)\Theta^-(\Lambda_k)\enddi and \begdi Q_k \triangleq \Theta^-(\Lambda_k)\Theta^+(\Lambda^\#_k) - \Theta^-(\Lambda^\#_k)\Theta^+(\Lambda_k). \enddi
\endth \vspace{0.1in}
\begpr The proof follows by direct application of Lemmas \ref{le:prop_thetas}, \ref{le:commutators_evolution_equation} and \ref{le:commutators_evaluated}. That is, equation \rref{eq:general_evolution_vector} is re-written using \rref{eq:commutator_xH}-\rref{eq:commutator_xLTL} as the following bilinear QSDE
\begeq \label{eq:physical_spin_evolution_evaluated}
\begin{split}
\lefteqn{ dx  =  - 2 \Theta^-(\alpha)x\,dt +\frac{4 \pmb{i}}{n} \sum_{k=1}^{n_w} \Theta^-(\Lambda^\#_k)\Lambda^T_k \,dt}  \\
      & + \sum_{k=1}^{n_w} \left( \Theta^-(\Lambda_k)\Theta^-(\Lambda^\#_k) + \Theta^-(\Lambda^\#_k)\Theta^-(\Lambda_k) \right)x\,dt \\
      & -\pmb{i} \sum_{k=1}^{n_w} \left( \Theta^-(\Lambda_k)\Theta^+(\Lambda^\#_k) - \Theta^-(\Lambda^\#_k)\Theta^+(\Lambda_k) \right) x\,dt \\
      & + \pmb{i}\left( \Theta^-(\Lambda_1^\#-\Lambda_1)x , \cdots , \Theta^-( \Lambda_{n_w}^\#-\Lambda_{n_w}) x\right)   \,d\bar{W}_1 \\
      & + \left( \Theta^-(\Lambda_1+\Lambda^\#_1)x ,  \cdots , \Theta^-(\Lambda_{n_w}+\Lambda^\#_{n_w}) x  \right) \,dW.
\end{split}
\endeq
Also, as mentioned in Section \ref{sec:section2}, the output fields $\bar{Y}_1$ and $\bar{Y}_2$ depend linearly on $L$, $L^\dagger$ and the input fields $\bar{W}_1$ and $\bar{W}_2$, i.e.,
\begdi
\left(\begin{array}{c}
d\bar{Y}_1 \\ d\bar{Y}_2
\end{array}\right)= \left(\begin{array}{c}
\Lambda + \Lambda^\# \\ \pmb{i}(\Lambda^\# - \Lambda) 
\end{array}\right) x \,dt + 
\left(\begin{array}{c}
d\bar{W}_1 \\ d\bar{W}_2
\end{array}\right). 
\enddi 
It is now easy to identify matrices $A_0,A,B_{1k}, B_{2k}, C_1$ and $C_2$, which ends the proof.
\endpr

Note that all matrices involved in the above equation are real. To confirm that, observe that $\Lambda^\#-\Lambda$ is purely imaginary and $\Lambda+\Lambda^\#$ is purely real. Now fix $k$ and compute the real part of $\Theta^-(\Lambda_k)\Lambda_k^\dagger$ and $\Theta^-(\Lambda_k)\Theta^+(\Lambda_k^\#) - \Theta^-(\Lambda_k^\#)\Theta^+(\Lambda_k)$. Given that $(\Theta^-(\Lambda_k)\Lambda_k^\dagger)^\#=-\Theta^-(\Lambda_k)\Lambda_k^\dagger$, one has that
\begdi Re\{\Theta^-(\Lambda)\Lambda^\dagger\}= \frac{1}{2}\left(\Theta^-(\Lambda)\Lambda^\dagger + (\Theta^-(\Lambda)\Lambda^\dagger)^\#  \right)=0.\enddi
Also, 
\begin{align*} 
\lefteqn{Re\{\Theta^-(\Lambda_k)\Theta^+(\Lambda_k^\#) - \Theta^-(\Lambda_k^\#)\Theta^+(\Lambda_k) \} } & \\
& = \frac{1}{2}\left(\Theta^-(\Lambda_k)\Theta^+(\Lambda_k^\#) - \Theta^-(\Lambda_k^\#)\Theta^+(\Lambda_k) \right. \\
 & \, \left. \rule{0in}{0.2in} \hspace*{0.1in} + \left(\Theta^-(\Lambda_k)\Theta^+(\Lambda_k^\#) - \Theta^-(\Lambda_k^\#)\Theta^+(\Lambda_k)\right)^\#  \right) =0.%\\
%& = \frac{1}{2}\left(\Theta^-(\Lambda_k)\Theta^+(\Lambda_k^\#) - \Theta^-(\Lambda_k^\#)\Theta^+(\Lambda_k) \right. \\
% & \, \left. \rule{0in}{0.2in} \hspace*{0.3in}+ \Theta^-(\Lambda_k^\#)\Theta^+(\Lambda_k) - \Theta^-(\Lambda_k)\Theta^+(\Lambda_k^\#)  \right) \\
%& = 0.
\end{align*}
% Therefore, $\pmb{i}(\Theta^-(\Lambda_k)\Theta^+(\Lambda_k^\#) - \Theta^-(\Lambda_k^\#)\Theta^+(\Lambda_k))$ and $\pmb{i}\Theta(\Lambda)\Lambda^\dagger$ are real valued. Since this argument is valid for all $k$, then all matrices in \rref{eq:physical_spin_evolution_evaluated} are real valued. Also, as mentioned in Section \ref{sec:section2}, the output fields $\bar{Y}_1$ and $\bar{Y}_2$ depend linearly on $L$, $L^\dagger$ and the input fields $\bar{W}_1$ and $\bar{W}_2$, i.e.,
% \begdi
% \left(\begin{array}{c}
% d\bar{Y}_1 \\ d\bar{Y}_2
% \end{array}\right)= \left(\begin{array}{c}
% \Lambda + \Lambda^\# \\ \Lambda^\# - \Lambda 
% \end{array}\right) x \,dt + \left(\begin{array}{cc}
% I_{n_w} & 0 \\ 0 & I_{n_w} 
% \end{array}\right) 
% \left(\begin{array}{c}
% d\bar{W}_1 \\ d\bar{W}_2
% \end{array}\right). 
% \enddi 
Moreover, by direct inspection $ {{B}_{ik}}^T = -{B}_{ik}$ for $i=1,2$ and $k=1,\hdots,n_w$. 

Now, from a control perspective, it is necessary to characterize when a bilinear QSDE posses underlying Hamiltonian and coupling operators which allows to express the matrices comprising \rref{eq:bilinear_system} and \rref{eq:bilinear_system_output} as in Theorem \ref{th:physical_realizability_def}. Thus, the second and most relevant result of the paper is given in the next theorem, which establishes necessary and sufficient conditions for the physical realizability of a bilinear QSDE.

\begth \label{th:physical_realizability}
System \rref{eq:bilinear_system} with output equation \rref{eq:bilinear_system_output} is physically realizable if and only if %it satisfies
\renewcommand*\theenumi{{\roman{enumi}}}
\renewcommand*\labelenumi{$\theenumi.$}
%\vspace*{0.05in}
\begin{enumerate}
\item \label{itm:theorem_physical_realizability1} $\displaystyle A_0=\frac{1}{n}\sum_{k=1}^{n_w}(\pmb{i}{B}_{1k} + {B}_{2k}) \left((C_1)_k+\pmb{i} (C_2)_k\right)^T$,
\vspace*{0.05in}
\item \label{itm:theorem_physical_realizability2} $\displaystyle {B}_{1k}= \Theta^-((C_2)_k)$,
\vspace*{0.1in}
\item \label{itm:theorem_physical_realizability21} $\displaystyle {B}_{2k}= \Theta^-((C_1)_k)$,
\vspace*{0.05in}
\item \label{itm:theorem_physical_realizability3} $\displaystyle A+A^T+\sum_{i,k=1}^{2,n_w} {B}_{ik} {{B}_{ik}}^T =\frac{n}{2}\Theta^+(A_0)$,
\vspace*{0.05in}
\end{enumerate}
where $(C_i)_k$ indicates the $k$-th row of $C_i$. In which case, the coupling matrix can be identified to be \begdi \Lambda=\frac{1}{2}(C_1+\pmb{i}C_2),\enddi
and $\alpha$, defining the system Hamiltonian, is
\begin{align} \label{eq:Hamiltonian_physical_realizability}
\nonumber  \lefteqn{\hspace*{-0.8in}\alpha = \frac{1}{4n}\vec\left(A^T-A +\frac{1}{2} \sum_{k=1}^{n_w}\left(  [B_{2k} , \Theta^+((C_2)_k)] \right. \right. } \\
 & \, \left. \left. \rule{0in}{0.2in} \hspace*{-0.6in} - [B_{1k}, \Theta^+((C_1)_k)] \right) \right)^T  F.
\end{align}
%where the commutator of matrices in \rref{eq:Hamiltonian_physical_realizability} is defined as usual, i.e., $[A,B]=AB-BA$ for $A,B\in \re^{s\times s}$.%the commutator of matrices 
\endth \vspace*{0.1in}

\begpr
Assuming that \rref{eq:bilinear_system} and \rref{eq:bilinear_system_output} are physically realizable implies that \rref{subeqn:F0}-\rref{subeqn:H2} are satisfied. By %straightforward 
comparison, conditions $\rref{itm:theorem_physical_realizability2}$-$\rref{itm:theorem_physical_realizability21}$ hold. Condition $\rref{itm:theorem_physical_realizability1}$ is written from \rref{subeqn:H12} and \rref{subeqn:H2} as
\begin{align*}
A_0 & = \frac{\pmb{i}}{n} \sum_{k=1}^{n_w}\Theta^-((C_1)_k-\pmb{i} (C_2)_k)((C_1)_k+\pmb{i} (C_2)_k)^T \\
    & = \frac{1}{n}\sum_{k=1}^{n_w}(\pmb{i}B_{1k} + B_{2k})((C_1)_k+\pmb{i} (C_2)_k)^T 
\end{align*}
Now, one has that
\begin{align*} %\label{eq:GGdagger_to_Lambda}
{B}_{1k} {{B}_{1k}}^T  = {} & \Theta^-(\Lambda_k^\#-\Lambda_k)^2 \\
= {} & \Theta^-(\Lambda^\#)\Theta^-(\Lambda^\#) - \Theta^-(\Lambda^\#)\Theta^-(\Lambda)   \\
 &   - \Theta^-(\Lambda) \Theta^-(\Lambda^\#) +\Theta^-(\Lambda)\Theta^-(\Lambda).
\end{align*}
Similarly,
\begin{align*} %\label{eq:GGdagger_to_Lambda_2}
B_{2k} B_{2k}^T  =  & -\Theta^-(\Lambda_k +\Lambda_k^\#)^2 \\
=  & -\Theta^-(\Lambda_k^\#) \Theta^-(\Lambda_k^\#) - \Theta^-(\Lambda_k^\#)\Theta^-(\Lambda_k) \\
 &   - \Theta^-(\Lambda_k) \Theta^-(\Lambda_k^\#) -\Theta^-(\Lambda_k)\Theta^-(\Lambda_k).
\end{align*}
Thus, ${B}_{1k} {{B}_{1k}}^T+{B}_{2k} {{B}_{2k}}^T  = -2 R_k $. 
%Define $Q_k=\Theta^-(\Lambda_k)\Theta^+(\Lambda^\#_k) - \Theta^-(\Lambda^\#_k)\Theta^+(\Lambda_k)$. 
One can now rewrite $A$ in terms of $\alpha,{B}_{1k}$ and ${B}_{2k}$ as 
\begin{align} \label{eq:F_physically_realizable}
\nonumber A = & -2\Theta^-(\alpha) \\
 &  - \frac{1}{2} \sum_{k=1}^{n_w}\left({B}_{1k}{{B}_{1k}}^T + {B}_{2k}{{B}_{2k}}^T\right)  -\pmb{i} Q_k.%\left( \Theta^-(\Lambda_k)\Theta^+(\Lambda^\#_k) - \Theta^-(\Lambda^\#_k)\Theta^+(\Lambda_k) \right)  .
\end{align}
Similarly, 
\begin{align} \label{eq:FT_physically_realizable}
\nonumber A^T = & \, 2\Theta^-(\alpha) \\
& -\frac{1}{2} \sum_{k=1}^{n_w} \left({B}_{1k}{{B}_{1k}}^T
+ {B}_{2k}{{B}_{2k}}^T\right) - \pmb{i} Q_k^T. 
\end{align}
%Adding \rref{eq:F_physically_realizable} and \rref{eq:FT_physically_realizable} 
%To obtain $\rref{itm:theorem_physical_realizability3}$, 
Adding \rref{eq:F_physically_realizable} and \rref{eq:FT_physically_realizable} gives %. It then follows that
\begdi%eq \label{eq:FT_physically_realizable}
A+A^T=- \sum_{k=1}^{n_w} \left({B}_{1k}{{B}_{1k}}^T + {B}_{2k}{{B}_{2k}}^T\right) - \pmb{i}\left( Q_k+Q_k^T \right). %\Theta^-(\Lambda_k)\Theta^+(\Lambda^\#_k) - \Theta^-(\Lambda^\#_k)\Theta^+(\Lambda_k) + \left( \Theta^-(\Lambda_k)\Theta^+(\Lambda^\#_k) - \Theta^-(\Lambda^\#_k)\Theta^+(\Lambda_k) \right)^T \right) 
\enddi  
The $(i,j)$ component of $Q_k+Q_k^T$ is computed as %obtained from %the last summands of \rref{eq:F_physically_realizable} and \rref{eq:FT_physically_realizable} is obtained from
\begin{align*}
\left(Q_k+Q_k^T\right)_{ij} %\left( \Theta^-(\Lambda_k) \Theta^+(\Lambda_k^{\#})-\Theta^-(\Lambda_k^\#) \Theta^+(\Lambda_k) - \Theta^+(\Lambda_k^\#) \Theta^-(\Lambda_k)+\Theta^+(\Lambda_k) \Theta^-(\Lambda_k^{\#})\right)_{ij} } \\
   = &\;  \Theta^-(\Lambda_k)\Theta^+(\Lambda^\#_k)  - \Theta^-(\Lambda^\#_k)\Theta^+(\Lambda_k)\\
&  - \Theta^+(\Lambda^\#_k)\Theta^-(\Lambda_k) + \Theta^+(\Lambda_k) \Theta^-(\Lambda^\#_k) \\
   = & -\Lambda_k^T F_iD_j \Lambda_k^\# + \Lambda_k^\dagger F_iD_j \Lambda_k \\
&  - \Lambda_k^\dagger D_iF_j \Lambda_k + \Lambda_k^T D_iF_j \Lambda_k^\#. 
\end{align*}
Note that every summand is a scalar, which is equal to its transpose. By \rref{eq:D_F_5} and \rref{eq:D_F_6}, it follows that  
\begin{align*} %\label{eq:antisymmetric_thetaplus_thetaminus}
\nonumber \lefteqn{\hspace*{-0.2in}\left(Q_k+Q_k^T\right)_{ij}} \\
%\nonumber & = \Lambda_k^T (D_iF_j - F_iD_j )\Lambda_k^\# + \Lambda_k^\dagger (D_iF_j-F_iD_j) \Lambda_k \\
\nonumber & = \Lambda_k^T (D_iF_j + D_jF_i )\Lambda_k^\# + \Lambda_k^T (F_iD_j + F_jD_i ) \Lambda_k^\# \\
%\nonumber & = \Lambda_k^T [D_i,F_j]\Lambda_k^\# - \Lambda_k^T [F_i,D_j] \Lambda_k^\# \\
\nonumber & = \Lambda_k^T \sum_{k=1}^{s}d_{jik}F_k\Lambda_k^\# + \Lambda_k^T \sum_{k=1}^{s}d_{ijk}F_k \Lambda_k^\#\\
\nonumber & = 2\sum_{k=1}^{s} d_{ijk}\Lambda_k^T D_k\Lambda_k^\# \\
& = 2 \left(\Theta^+(\Theta^-(\Lambda_k)\Lambda^\#_k )\right)_{ij}\\
& = -2 \left(\Theta^+(\Theta^-(\Lambda^\#_k)\Lambda_k )\right)_{ij}\\
& = \frac{n\pmb{i}}{2}\left(\Theta^+(A_0)\right)_{ij}. 
\end{align*}
%which implies that the matrix $Q_k=\Theta^-(\Lambda_k) \Theta^+(\Lambda_k^{\#})-\Theta^-(\Lambda_k^\#) \Theta^+(\Lambda_k)$ is antisymmetric. Observe that this antisymmetry is independent of the vectors applied to $\Theta^-$ and $\Theta^+$. 
Therefore, adding \rref{eq:F_physically_realizable} and \rref{eq:FT_physically_realizable} gives
\begin{align*}
 A+A^T + \sum_{i,k=1}^{2,n_w}  B_{ik} B_{ik}^T = \frac{n}{2}\Theta^+(A_0),
\end{align*}
which is condition $\rref{itm:theorem_physical_realizability3}$. Conversely, one needs to show that if conditions $\rref{itm:theorem_physical_realizability1}$-$\rref{itm:theorem_physical_realizability3}$ of Theorem \ref{th:physical_realizability} are satisfied, then there exist matrices $\alpha$ and
$\Lambda$ such that system \rref{eq:bilinear_system} is physically realizable. Let 
\begin{align} \label{eq:realization_hamiltonian}
\nonumber \Theta^-(\alpha) \triangleq & \,  \frac{1}{4} \left( \rule{0in}{0.25in} A^T - A \right. \\
& \hspace*{-0.5in} \left. + \frac{1}{2}\sum_{k=1}^{n_w} \left( [B_{2k} , \Theta^+((C_2)_k)] -  [B_{1k} ,\Theta^+((C_1)_k)] \right) \right).
\end{align}
It is trivial to check that the right-hand-side of \rref{eq:realization_hamiltonian} is antisymmetric and hence this equation uniquely defines $\alpha$ via \rref{eq:Theta_definition}. Also, let \begeq \label{eq:def_lambda_Theorem}\Lambda=\frac{1}{2}(C_1+\pmb{i}C_2).\endeq Then, $Q_k$ can be written in terms of $B_1,B_2,C_1,$ and $C_2$ as follows 
% \begin{align*}
%       {G}_{1k} &= \Theta^-\left({\pmb i}(\Lambda^\#_k - \Lambda_k)\right), \\ 
%       {G}_{2k} &= \Theta^-(\Lambda_k + \Lambda^\#_k),  \\
%       H_1   &=  \Lambda+\Lambda^\#,  \\
%       H_2   &=  \pmb{i} \left( \Lambda^\#-\Lambda \right). 
% \end{align*}
\begin{align*}
%       \lefteqn{ \hspace*{-0.2in}\Theta^-(\Lambda_k) \Theta^+(\Lambda_k^{\#})-\Theta^-(\Lambda_k^\#) \Theta^+(\Lambda_k) } \\
Q_k = & \, \frac{1}{4}\left( \Theta^-((C_1+\pmb{i}C_2)_k) \Theta^+((C_1-\pmb{i}C_2)_k) \right.  \\
& \,\left. -\Theta^-((C_1-\pmb{i}C_2)_k) \Theta^+((C_1+\pmb{i}C_2)_k) \right)\\
& \, - \Theta^+((C_1+\pmb{i}C_2)_k) \Theta^-((C_1-\pmb{i}C_2)_k)\\
& \, +\Theta^+((C_1-\pmb{i}C_2)_k) \Theta^-((C_1+\pmb{i}C_2)_k)\\
%  = & \, -\frac{\pmb{i}}{4}\left( \Theta^-((C_1)_k) \Theta^+((C_2)_k) - \Theta^-((C_2)_k) \Theta^+((C_1)_k) \right. \\
%    & \, \left.  + \Theta^-((C_1)_k) \Theta^+((C_2)_k) - \Theta^-((C_2)_k) \Theta^+((C_1)_k)  \right) \\
= & \, -\frac{\pmb{i}}{2}\left( \Theta^-((C_1)_k) \Theta^+((C_2)_k) \right. \\
& \, \left. \hspace*{0.34in} - \Theta^-((C_2)_k) \Theta^+((C_1)_k) \right). % \\
%= & \, -\frac{\pmb{i}}{2}\left( B_{2k}  \Theta^+((C_2)_k) - B_{1k} \Theta^+((C_1)_k) \right).
\end{align*}
From $\rref{itm:theorem_physical_realizability2}$ and $\rref{itm:theorem_physical_realizability21}$, it follows that
\begin{align*}
Q_k & =\Theta^-(\Lambda_k)\Theta^-(\Lambda^\#_k) + \Theta^-(\Lambda^\#_k)\Theta^-(\Lambda_k) \\
% -\frac{\pmb{i}}{2}\left( B_{2k}  \Theta^+((C_2)_k) - B_{1k} \Theta^+((C_1)_k) \right).% \\
& =  -\frac{\pmb{i}}{2}\left( B_{2k}  \Theta^+((C_2)_k) - B_{1k} \Theta^+((C_1)_k) \right).
\end{align*}
Then.
\begin{align*}
Q^T_k - Q_k & = \frac{\pmb{i}}{2}\left( [B_{2k} , \Theta^+((C_2)_k)] - [B_{1k} ,\Theta^+((C_1)_k)] \right).% \\
%& =  -\frac{\pmb{i}}{2}\left( B_{2k}  \Theta^+((C_2)_k) - B_{1k} \Theta^+((C_1)_k) \right).
\end{align*}
Similarly, it is simple to write $R_k$ in terms of $C_1$ and $C_2$. That is, %This symmetry is only consequence of mappings $\Theta^-$ and $\Theta^+$ and not dependent on their arguments. Writing $R_k$ in terms of $C_1$ and $C_2$ gives %Therefore, it is sufficient to compute 
\begin{align*}
%\lefteqn{\hspace*{-0.2in}\Theta^-(\Lambda_k)\Theta^-(\Lambda^\#_k) + \Theta^-(\Lambda^\#_k)\Theta^-(\Lambda_k)} & \\
R_k = &  \, \frac{1}{4} \left( \Theta^-((C_1+\pmb{i}C_2)_k)\Theta^-((C_1-\pmb{i}C_2)_k) \right. \\
  & \,\left. \hspace*{0.15in} + \Theta^-((C_1-\pmb{i}C_2)_k)\Theta^-((C_1+\pmb{i}C_2)_k) \right)\\
 = &  \, \frac{1}{4}\left( \Theta^-((C_1)_k) \Theta^-((C_1)_k) - \pmb{i}\Theta^-((C_1)_k) \Theta^-((C_2)_k) \right. \\
   &  \, + \,\pmb{i}\Theta^-((C_2)_k) \Theta^-((C_1)_k)+\Theta^-((C_2)_k) \Theta^-((C_2)_k) \\
   &  \, + \Theta^-((C_1)_k) \Theta^-((C_1)_k) + \pmb{i}\Theta^-((C_1)_k) \Theta^-((C_2)_k)  \\
   &  \, \left. - \,\pmb{i}\Theta^-((C_2)_k) \Theta^-((C_1)_k)+\Theta^-((C_2)_k) \Theta^-((C_2)_k) \right) \\
= &  \, \frac{1}{2}\left( \Theta^-((C_1)_k) \Theta^-((C_1)_k)  + \Theta^-((C_2)_k) \Theta^-((C_2)_k) \right).
%= &  \, -\frac{1}{2}\left( G_{1k} G_{1k}^T  +  G_{2k} G_{2k}^T \right).
\end{align*}
% Recall from Section \ref{sec:section2} that $H_i={H_i}^\#$ and $G_i={G_i}^\#$ for $i=1,2$. It then follows that  
% \begin{align*}
% \lefteqn{\Lambda^\dagger \Lambda+\Lambda^T \Lambda^\# - 2\Lambda \Lambda^\dagger I} \\
% & = \frac{1}{2}\left( H_1^\dagger H_1 +H_2^\dagger H_2 - \left(H_1 H_1^\dagger + H_2 H_2^\dagger \right)I  \right).
% \end{align*}
It clear that $R_k$ is symmetric. From $\rref{itm:theorem_physical_realizability2}$ and $\rref{itm:theorem_physical_realizability21}$, one obtains that 
\vspace*{-0.05in}\begin{align*}
R_k & =\Theta^-(\Lambda_k)\Theta^-(\Lambda^\#_k) + \Theta^-(\Lambda^\#_k)\Theta^-(\Lambda_k) \\
& = -\frac{1}{2}\left( B_{1k} B_{1k}^T  +  B_{2k} B_{2k}^T \right).
\end{align*}
From $\rref{itm:theorem_physical_realizability1}$ and \rref{eq:def_lambda_Theorem},
\begin{align*}
\Theta^+(A_0) & = -\frac{2\pmb{i}}{n} \left(Q_k +Q_k^T\right)  
\end{align*}
Since $\rref{itm:theorem_physical_realizability3}$ implies 
\begdi 
A^T= - A - \sum_{i,k=1}^{2,n_w} {B}_{ik} {{B}_{ik}}^T  + \frac{n}{2}\Theta^+(A_0),
\enddi 
one can use \rref{eq:realization_hamiltonian} to obtain
\begin{align*}
\lefteqn{\Theta^-(\alpha)} \\
 & =  \frac{1}{4}\left(-2A + \frac{n}{2}\Theta^+(A_0) - \sum_{k=1}^{n_w}  \left(\rule{0in}{0.18in} {B}_{1k} {{B}_{1k}}^T + {B}_{2k} {{B}_{2k}}^T \right. \right. \\
& \, \left. \rule{0in}{0.25in} \left. \hspace*{0.15in}  + \frac{1}{2}\left([B_{2k},  \Theta^+((C_2)_k)] - [B_{1k} , \Theta^+((C_1)_k)]\right) \right) \right)\\ 
 & =	 -\frac{1}{2} A - \frac{1}{4}\sum_{k=1}^{n_w} \left( 2 R_k - \pmb{i} (Q_k -Q_k^T) - \pmb{i} (Q_k^T -Q_k) \right)\\
 & =	 -\frac{1}{2} A - \frac{1}{2}\sum_{k=1}^{n_w} \left(R_k - \pmb{i} Q_k \right),
\end{align*}
which is equivalent to \rref{subeqn:F}. Moreover, using \rref{eq:E^T_by_E}, \rref{eq:realization_hamiltonian} and applying the stacking operator to $\Theta^-(\alpha)$, $\alpha$ is explicitly obtained as %applying the stacking operator and multiplying by $E^T$, which gives % operator and $(\Theta(\alpha))=E\alpha^T=\frac{1}{4}\vec\left(A^T-A - \frac{1}{2}\left( B_{2k}  \Theta^+((C_2)_k) - B_{1k} \Theta^+((C_1)_k) \right)\right)$. Multiplying both sides by $E^T$ leaves 
\begdi
\left(F^T\vec\left( \Theta^-(\alpha)\right) \right)^T = \left(F^T F\alpha^T \right)^T = n\alpha.
\enddi
Hence,\vspace*{-0.05in}
\begin{align*}
\lefteqn{\alpha= \frac{1}{4n} \vec\left(A^T-A + \frac{1}{2}\sum_{k=1}^{n_w} [B_{2k} , \Theta^+((C_2)_k)]     \right. } \\
 & \, \left.  \rule{0in}{0.22in} \hspace*{0.7in} - [B_{1k}, \Theta^+((C_1)_k)]  \right)^T  F,
\end{align*}
which completes the proof. 
\endpr
\vspace*{0.05in}

\noindent Note that the physical realizability conditions do not require the computation of the Hamiltonian \rref{eq:Hamiltonian_physical_realizability},  which depends on the structure constants $d$ and $f$, in order to know whether or not the system given by equations \rref{eq:bilinear_system} and \rref{eq:bilinear_system_output} is quantum.

\section{Conclusions} \label{sec:conclusions}

A condition for physical realizability was given for open multi-level quantum systems. Under this condition it was shown that there exist operators $\mathcal{H}$ and $L$ such that the bilinear QSDE \rref{eq:bilinear_system} with output equation \rref{eq:bilinear_system_output} can be written as in \rref{eq:general_evolution_vector}. This condition used explicitly the algebra generated by $SU(n)$. Moreover, the interaction of the system with multiple quantum fields was introduced to the formalism. %The result was presented as a generalization of the condition for two-level systems  results given  

\end{document}